# PorePy: An Open-Source Software for Simulation of Multiphysics Processes in Fractured Porous Media

Eirik Keilegavlen[1], Runar Berge[1], Alessio Fumagalli[1], Michele Starnoni[1], Ivar Stefansson[1], Jhabriel Varela[1], Inga Berre[1]

[1]Department of Mathematics, University of Bergen.

Corresponding author: Eirik.Keilegavlen@uib.no, Department of Mathematics, University of Bergen, Pb 7800, 5020 Bergen, Norway.

## Abstract

Development of models and dedicated numerical methods for dynamics in fractured rocks is an active research field, with research moving towards increasingly advanced process couplings and complex fracture networks. The inclusion of coupled processes in simulation models is challenged by the high aspect ratio of the fractures, the complex geometry of fracture networks and the crucial impact of processes that completely change characteristics on the fracture-rock interface. This paper provides a general discussion of design principles for introducing fractures in simulators, and defines a framework for integrated modeling, discretization and computer implementation. The framework is implemented in the simulation software PorePy, which can serve as a flexible prototyping tool or multiphysics problems in fractured rocks. Based on a representation of the fractures and their intersections as lower-dimensional objects, we discuss data structures for mixed-dimensional meshes, formulation of multiphysics problems and discretizations that utilize existing software. We further present the implementation of these concepts in the PorePy open-source software tool, which is aimed at coupled simulation of flow and transport in three-dimensional fractured reservoirs as well as deformation of fractures and the reservoir in general. We present validation by benchmarks for flow, poroelasticity and fracture deformation in fractured porous media. The flexibility of the framework is then illustrated by simulations of fully coupled flow and transport and of injection driven deformation of fractures. All results reported herein can be reproduced by openly available simulation scripts.

## 1. Introduction

Simulation of flow, transport and deformation of fractured rocks is of critical importance to several applications such as subsurface energy extraction and storage and waste disposal. While the topics have received considerable attention the last decade, the development of reliable simulation tools remains a formidable challenge. Many reasons can be given for this, we here pinpoint four possible causes: First, while natural fractures are thin compared to the characteristic length of the domains of interest, their extent can span through the domain of interest [1]. The high aspect ratios make the geometric representation of fractures in the simulation model challenging. Second, the strongly heterogeneous properties of fractures compared to the matrix with respect to flow and mechanics call for methods that can handle strong parameter discontinuities as well as different governing physics for the fractures and the matrix, *e.g.* [2]–[4]. Third, phenomena of practical interest tend to involve multiphysics couplings, such as interaction between flow, temperature evolution, geo-chemical effects and fracture deformation [5]. Correspondingly, there is an ongoing effort to develop and introduce multiphysics couplings within simulation models [6]. Fourth, fracture networks have highly complex intersection geometries, which must be accounted for in the simulation models. We emphasize that, although the geometry of the walls of individual fractures can be complex by themselves, we will not consider this in any detail, but rather assume that averaged apertures etc. are available at the scale of discretizations.

Traditionally, simulation of flow-driven dynamics in fractured media has been based on two conceptual models: First, in an upscaled representation, the fracture network geometry and dynamical processes taking place in the network are replaced by equivalent continuum models, which resemble those used in non-fractured porous media. As these models do not resolve the fracture geometry, they are computationally efficient, and have been extended to cover a wide range of multiphysics couplings, as exemplified by the TOUGH2 family of codes [7], PFLOTRAN [8], and also *e.g.* [9]. The accuracy of the simulation is however highly dependent on the quality of the upscaled model, which in turn depends on the fractured domain's resemblance of a continuous medium with respect to the nature of the physical processes. In practice, the upscaling process ranges from treatable by analytical means for simple fracture geometries and dynamics [10], [11], to extremely challenging in the case of multiphysics couplings and complex fracture geometries [12], [13].

The second traditional class of models, known as the discrete fracture network (DFN) models, is constructed using an explicit representation of the fracture network in the simulation model, while ignoring the surrounding rock mass. The models combine highly accurate representation of dynamics in the fractures with computational efficiency from not having to deal with the rock matrix, which is highly desirable *e.g.* for fast model evaluation. DFN simulation models with a high level of sophistication have been developed, notably for coupled flow and transport, see for instance [14]–[16]. By themselves, DFN models cannot represent processes outside the fracture network; however, the models can be combined with continuum models to achieve fracture-matrix couplings.

The respective limitations of upscaled and DFN approaches have over the last decade led to an increased interest in the class of discrete fracture matrix (DFM) models. In DFM models, the fractures are sorted in two classes according to their importance for the dynamics in question [17]. The most important fractures are represented explicitly, while upscaled models are applied for the remaining fractures and the host rock. As such, DFM models represent a flexible compromise between upscaling and explicit representations. The models can represent governing equations in the rock matrix, fractures, and generally also in the intersections between fractures. For computational efficiency, it is common to represent fractures and their intersections as lower-dimensional objects embedded in the three-dimensional rock matrix [18], [19]. We refer to this as a mixed-dimensional model [20], and conversely refer to a model of a domain where only a single dimension is considered as fixed-dimensional.

DFM models can further be divided into two subgroups, according to whether they explicitly represent the fracture surfaces in the computational grid [17]. Models that apply non-conforming meshing include the Embedded Discrete Fracture Matrix model (EDFM) [21], and extended finite element methods (XFEM) [22], [23]. These methods avoid the complexities of conforming mesh generation discussed below, but must instead incorporate the fracture-matrix interaction in what becomes complex modifications of the numerical method for XFEM [24], or by constructing upscaled representation reminiscent of the challenge in continuum-type models [25]. For this reason, our interest herein is DFM methods with conforming meshes. By now, this type of DFM models have been developed for flow and transport, as well as mechanics and poroelasticity. Simulation models that incorporate DFM principles include DuMuX [26], CSMP [27], MOOSE-FALCON [28], [29], OpenGeoSys [30] and Flow123d [31].

The utility of a rapid prototyping framework is illustrated by the wide usage of the Matlab Reservoir Simulation Toolbox (MRST) [32], [33], mainly for non-fractured porous media. Similarly, research into strongly coupled processes in mixed-dimensional geometries will benefit from software of similar flexibility and with a structure tailored to the specific challenges related to fractured porous media.

The goal of this paper is two-fold: First, we review challenges related to design of simulation frameworks for multiphysics couplings in mixed-dimensional geometries. Our aim is to discuss design choices that must be made in the implementation of any DFM simulator, including data structures for mixed-dimensional geometries, and representation and discretization of multiphysics problems. Second, we describe a framework for integrated modeling, discretization and implementation, and an open-source software termed PorePy adhering to this framework. Key to our approach is a decomposition of the geometry into separate objects for rock matrix, individual fractures and fracture intersections. Governing equations can then be defined separately on each geometric object, as well as on the connection between the objects. This allows for significant code reuse from the discretization of fixed-dimensional problems; thus, our design principles are also applicable to more general PDE software frameworks, such as FEniCS [34], Dune [35] and FireDrake [36]. Furthermore, for scalar and vector elliptic problems (flow and deformation) the models rest on a solid mathematical formulation [37]–[39].

Built on the object-based mixed-dimensional geometry, PorePy offers several discretization schemes for mathematical models of common processes, such as flow, transport and mechanical deformation. Multiphysics couplings are easily formulated, and their discretization depends on the availability of appropriate discretization schemes. Moreover, the framework allows for different geometric objects to have different primary variables and governing equations. The software can be used for linear and non-linear problems, with the latter treated by automatic differentiation. For DFM models that explicitly represent the fractures in the computational grid, meshing is a major technical challenge, in particular for 3d problems, PorePy offers automatic meshing of fractured domains in 2d and 3d, relying on the third-party software Gmsh to construct the mesh [40]. The software is fully open-source (see [www.github.com/pmgbergen/porepy](www.github.com/pmgbergen/porepy)) and is released under the GNU General Public License (GPL).

The paper is structured as follows: In Section 2, we present the principles whereupon we have built the mixed-dimensional framework in PorePy. Section 3 deals with modeling and discretization of physical processes central to fractured porous media: single-phase flow, heat transport, poroelastic rock deformation, and fracture deformation modeled by contact mechanics. In Section 4, we benchmark our approach and the library PorePy by well-established test cases. In Section 5, we present two complex examples to illustrate the potential of the framework with respect to advanced physical processes, followed by conclusions in Section 6.

# 2 Design principles for mixed-dimensional simulation tools

A simulation model for a specific dynamical process in mixed-dimensional media requires three main ingredients: A representation of the mixed-dimensional geometry, governing equations for dynamics within and between the geometric objects (rock matrix, fractures, their intersections), and a strategy for discretization and assembly of the equations on the geometry. On the more fundamental level of simulator design, important questions to clarify include how much of the mixed-dimensional geometry to include, which type of couplings between different geometric objects to permit, and how to establish communication between the geometric objects.

In this section, we discuss principles for process couplings in a general context of fractured rocks, together with representation of the geometry in a continuous and discrete setting. As we will see, the design choices cannot be done independently, for instance the coupling structure puts constraints on the representation and data structure for the geometry. We further present the specific models underlying PorePy, including mesh generation, construction of projection operators between geometric entities, and discretization and assembly on mixed-dimensional geometries. The general discussion herein is supplemented by concrete examples of modeling and discretization of important governing processes presented in Section 3.

## 2.1 Representation of a mixed-dimensional geometry

We consider the geometry of a fracture network embedded in a 3d domain; 2d domains are treated by the natural simplification. In general, the geometry formed by the fracture network consists of objects of dimension 2 (the fractures), 1 (fracture intersections) and 0 (intersections of intersection lines), in addition to the 3d domain itself. An important decision for the modeling of dynamics in the domain is which parts of the geometry to represent in the model. We emphasize that as our focus herein is DFM models, it is assumed that at least the fractures in question will be explicitly represented in the simulation model, and furthermore that the simulation grid will conform to the fractures.

We differ between two approaches to representation of the fracture geometry: The first explicitly represents the full hierarchy of geometric objects (3d-0d) as described above. However, for many processes, one can to a good approximation assume that the main dynamics take place in the matrix or in the fractures, while objects of co-dimension more than 1 (intersection lines and points) mainly act as transition zones between fractures. This observation motivates the other approach: The matrix and fractures are represented explicitly, together with some model for direct fracture-fracture interaction.

Representation only of matrix and fractures and not the intersections in some sense constitutes the minimal modification to an existing fixed-dimensional model and has been a popular choice *e.g.* for flow and transport problems [41]. The strategy has also been taken a long way towards practical applications, see for instance [42]. There are however drawbacks, notably in the treatment of fracture intersections: Without explicit access to the intersection objects, modeling of interaction between two fractures can be challenging. Significantly, the difficulties tend to increase with increasing complexity of the dynamics, such as countercurrent flow due to gravity and capillary forces, and when transitioning from 2d domains to 3d (i.e. the dimension of the intersections increases from zero to one). This has important consequences for model and method development, as issues related to *ad hoc* treatment of intersection dynamics may not manifest until relatively late in the development process.

Our preferred solution is to apply an equal representation of all geometric objects, independent of their dimension. This allows for flexible modeling of dynamics within all objects, and as we will see below, the implementation of couplings between geometric objects can be made independent of the objects' dimensions. The design choice has further advantages in terms of reuse of discretizations, as will be discussed in Section 3. In this framework, variables follow the domain decomposition approach and are associated with single subdomains or interfaces. This is followed through in our implementation, where the solution vectors in different subdomains are represented by different objects, even if they represent parts of the same physical quantity, for instance pressure or temperature.

Our approach to the geometry representation is illustrated in Figure 1, which shows the decomposition of a mixed-dimensional geometry into a hierarchy of geometric objects with accompanying meshes.

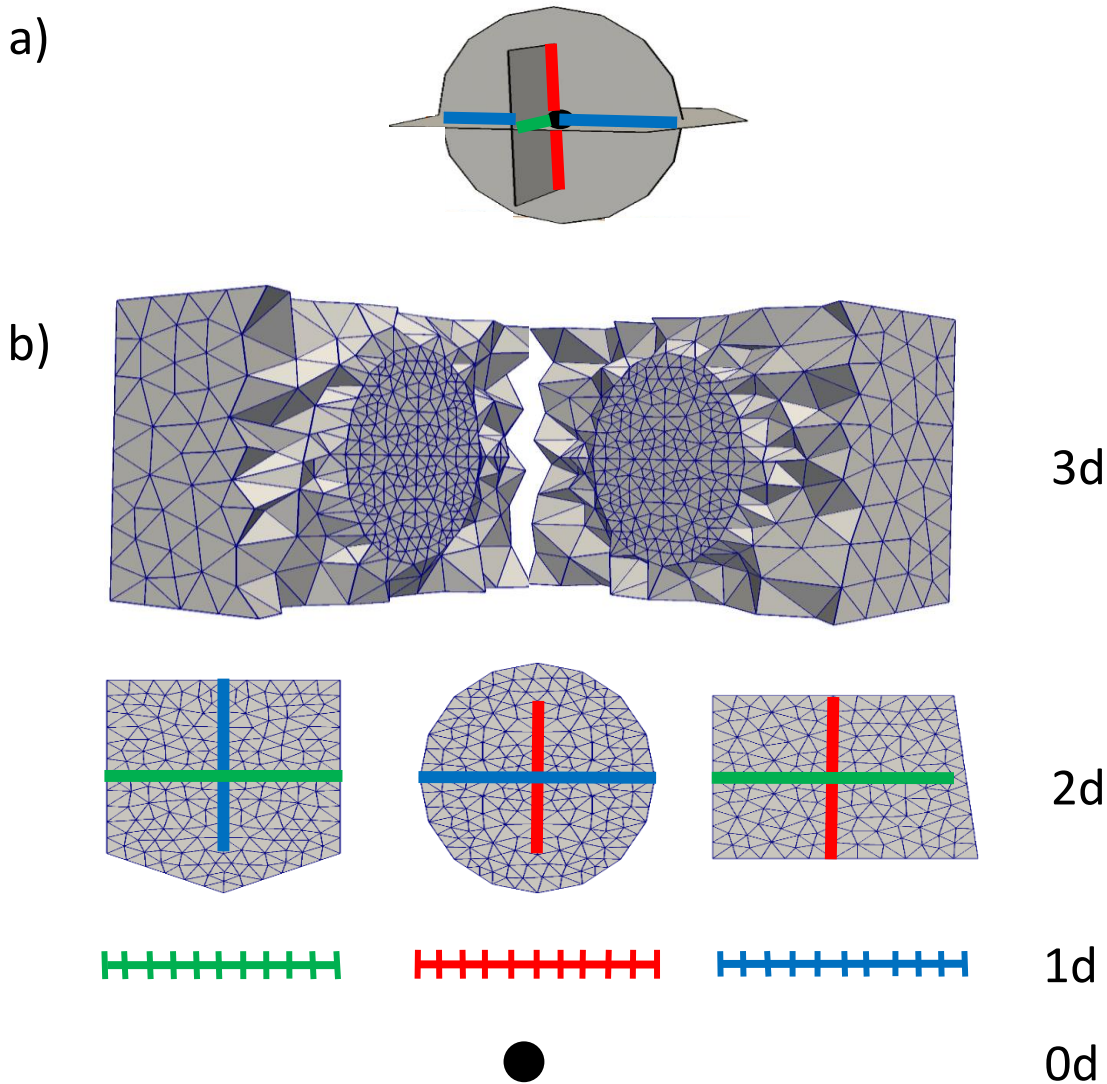

*Figure 1: Conceptual illustration of a fracture network, including meshing and lower-dimensional representation. a) Fracture network, the rock matrix is not visualized. b) Meshes of all subdomains. Fracture intersections (1d) are represented by colored lines, the 0d grid by a black circle. The 3d mesh is cut to expose the circular fracture.*

## 2.2 Permissible coupling structures between geometric objects

For modeling purposes, and for the design of data structures discussed next, it is important to establish which types of couplings between subdomains are permitted. In our framework we impose the following constraints on the modeling:

1. Only coupling between subdomains that are exactly one dimension apart is allowed.
2. Interaction between subdomains should be formulated as models associated with the interface between the subdomains. This requires interface variables that represent the interaction between the subdomains.
3. An interface law can depend on variables on the interface and the immediate subdomain neighbors, but not on variables associated with other subdomains or interfaces.

These choices have several important consequences: First, our framework explicitly rules out direct 3d-1d couplings. Although some of the ingredients presented herein could be of interest for such high dimensional gaps, notably a modified version of the mixed-dimensional grid structure presented below, the mathematical structure of the two problems is significantly different, and we have therefore not pursued a unified treatment. Second, our model does not permit direct coupling between objects of the same dimension, say, two fractures; the communication must go via a lower- or higher-dimensional object. Third, the restriction that couplings should be formulated in terms of interface variables makes the structure of the equations on a subdomain relatively simple, as the dynamics depends only on variables internal to the subdomain and on neighboring interfaces.

## 2.3 Data structure for geometry and data

The restrictions put on the coupling structure between subdomains give important guidance on the definition of a data structure for the compound mixed-dimensional grid, consisting of all subdomains and interfaces. As subdomains should only see neighboring interfaces, and the interfaces only communicate with their two neighboring subdomains, the grid can be represented in the simulation model as a graph, with the subdomains forming nodes, while interfaces are edges.

The graph is now the natural place to store all kind of data relating to the simulation including grids, parameters and variables, as well as information on which equations to solve in each subdomain, and which discretization schemes to apply. The data structure is thus a natural framework for defining advanced simulation models, examples will be shown in Sections 4 and 5. While this flexibility can be gained by applying domain decomposition to any problem [43], for fractured domains it comes as an added feature from what is already a natural data structure for the geometry.

## 2.4 Meshing and projections

Having defined the data structure for the mixed-dimensional geometry, we proceed to discussing meshing of the geometric objects and establishing projection operators for communication between the objects.

### 2.4.1 Mesh construction in mixed-dimensional geometries

A major technical difficulty of conforming DFM models is the construction of meshes. Obtaining meshes that conform to all geometric objects requires first, identification of all intersection lines and points, then meshing of objects of all dimensions, and finally identification of neighboring cells and faces on different domains, so that inter-object interaction can be modeled. In principle, the computation of fracture intersections is straightforward, following for instance [44]. In practice, this requires (automatic) decisions on when two objects should be considered distinct in the computational mesh; for complex networks this can be rather challenging. Notably, the question of whether objects should be considered spatially separated must be seen in connection with the prescribed mesh size, which puts practical constraints on how fine details can be resolved.

From a geometric description with all intersections identified, meshes of all objects can be constructed; in PorePy this is handled by a backport to Gmsh [40]. As post processing of the Gmsh output, standard simulation meshes are generated for all subdomains: The mesh for the 3d subdomain consists of all 3d cells, while for each fracture a mesh is composed of all faces of the 3d grid that lie on the surface of the fracture. Similarly, meshes for 1d intersection lines are formed by edges of the 3d grid that coincide with the line, while point meshes for 0d intersections are identified by nodes in the 3d grid. For 2d domains with 1d fractures and 0d intersections, the construction is similar. We emphasize that each of the meshes is implemented as a standard fixed-dimensional mesh, so that when a discretization scheme is applied to a subdomain, this is indistinguishable from the traditional fixed-dimensional operation. In this spirit, the grid structure used for individual meshes is agnostic to spatial dimension, with an implementation heavily inspired by that of MRST [33]. This grid structure in many cases facilitates an implementation of discretization schemes which is independent of dimension.

The meshes generated by Gmsh match between the subdomains. Moreover, the PorePy interface to Gmsh is restricted to simplex cells as these are most relevant for complex geometries. Non-matching grids can be introduced to PorePy by replacing meshes on individual subdomains; examples of computations on non-matching grids are given in Sections 4.1 and 5.1.

We finally note that tuning of mesh sizes in parts of the domain so that the resulting grid both resolves the local geometry and provides the desired accuracy for numerical computations can be a delicate task. Within PorePy, we attempt to handle this by setting a minimal mesh size and target mesh sizes for the fractures and the global boundary (far-field conditions). Based on these three user-provided entries, mesh size parameters are computed for all points in the fracture geometry and provided as guidance to Gmsh. In practice, Gmsh may override the settings, but nevertheless, the mesh size tuning combined with the automatic processing of the fracture geometry is a major capability of PorePy.

### 2.4.2 Mortar grids and projection operators

In addition to meshing on the subdomains, the interfaces are assigned separate meshes. These are used for discretization of the interface variables and serve as mortar grids for the coupling between the subdomains. Specifically, as the mortar grids allow for non-matching grids between subdomains, computational speedups can be achieved by combining fine grids in fractures, which are often the main venue for dynamical processes, with relatively coarse grids in the matrix. When using mortar technology to combine non-matching grids, non-uniform discretizations or physics, it is important to carefully design the mortar space so that the coupling does not introduce instabilities, see *e.g.* [45].

Transfer of variables between an interface and its neighboring subdomain is handled by projection operators. In the subsequent parts we will apply four different classes of projections. We have the mapping from an interface to the related subdomains indicated with $\Xi$, with a subscript indicating the index of the interface and a superscript the index of the subdomain, see Figure 2 for an illustration. We also introduce the projection operators from neighboring subdomains of an interface to the interface itself, denoted by the symbol Π with the same convention as before for sub- and superscripts. The actual definition of these objects is scope dependent and it will be specified when needed. In our implementation, we have only considered projections of lowest order which can be constructed by identifying overlapping areas between cells in the interface grid and faces and cells in the neighboring meshes. The actual construction of the projection needs to consider the nature of the variable to project, being of intensive or extensive kind.

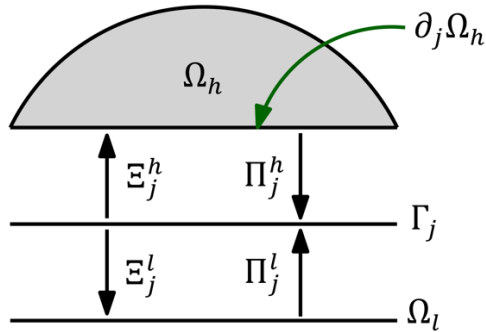

Figure 2: Generic geometry of a coupling between subdomains: An interface $\Gamma_j$ is coupled to a higher-dimensional subdomain $\Omega_h$ and a lower-dimensional subdomain $\Omega_l$. The projection operators between interfaces and subdomains are denoted by $\Pi$ (interface to subdomain) and $\Xi$ (subdomain to interface) with subscripts indicating the interface and superscript indicating the subdomain. In practice, $\Gamma_j$ will coincide with an internal boundary of $\Omega_h$, which we will refer to as $\partial_j \Omega_h$.

## 2.5 Governing equations and discretization

With the above framework, the task of defining governing equations and their discretization is split into two operations that to a large degree are independent, although this naturally depends on the physical process to be modeled. First, on the subdomains, the governing equations can often be defined (and to a large degree discretized) as if the problem were fixed-dimensional, while interaction with the interface variables takes the form of boundary conditions and source terms or body forces. These are terms that can be handled by any standard numerical method, and so the coupling structure paves the way for considerable reuse of existing simulation code designed for fixed-dimensional problems. The second operation involves coupling conditions on the interface, including projections of variables on the neighboring subdomains. This operation generally has no clear parallels for fixed-dimensional problems.

The details of the discretization can vary substantially depending on the governing equations and designated discretization schemes. We will give several examples of this in Section 3.

## 2.6 Global assembly

A global system of equations can be assembled from the components on individual subdomains and interfaces. The form and proper treatment of these equations differ according to whether the problem is stationary or time-dependent, linear or non-linear, but some ingredients of the implementation and structure of the problem are common. Specifically, for multiphysics problems with more than one primary variable, the global system of equations has a double block structure: One set of blocks stems from the geometric division into subdomains and interfaces. Within each subdomain and interface, there is a second set of blocks, with one block per variable. Access to this information is useful for design of tailored preconditioners and linear solvers, as well as post processing and visualization. PorePy has implemented a global degree-of-freedom manager, that keeps track of the block structure of the system of equations, as well as the numbering of individual degrees of freedom.

For visualization, an export filter to Paraview [46] is available. To aid analysis of simulation results, the export preserves the link between the data and its associated dimensions.

# 3 Modeling, discretization and implementation

In this section we apply the general framework presented above to three sets of governing equations, each of which is of high relevance for fractured porous media: The elliptic pressure equation, fully coupled flow and transport, and inelastic deformation of fractures due to poroelastic effects. As most of these processes are well established for fixed-dimensional, partly also for mixed-dimensional, problems, our main purpose is to cast the methodology in the general framework of Section 2, with discussions of modeling, extension of numerical methods designed for fixed-dimensional problems, and implementation aspects. The presentation will emphasize these three ingredients in varying degrees, with the aim that the section in total should illustrate the full power of the modeling framework.

The notation used for variables and subdomains is fixed as follows: Let $\Omega_i$ denote a generic subdomain, with variables in $\Omega_i$ marked by subscript $i$. A generic interface between two subdomains is represented by $\Gamma_j$, with subscript $j$ identifying interface variables. For a subdomain $\Omega_i$, the set of neighboring interfaces is split into interfaces towards subdomains of higher dimensions, denoted $\hat{S}_i$, and interfaces towards subdomains of lower dimensions, represented by $\check{S}_i$, see Figure 3. When discussing subdomain couplings for an interface $\Gamma_j$, we let the higher- and lower-dimensional subdomain be represented by $\Omega_h$ and $\Omega_l$, respectively, and associate variable subscripts $h$ and $l$, see Figure 2. Finally, the part of the boundary of $\Omega_h$ that coincides with $\Gamma_j$ is denoted $\partial_j \Omega_h$.

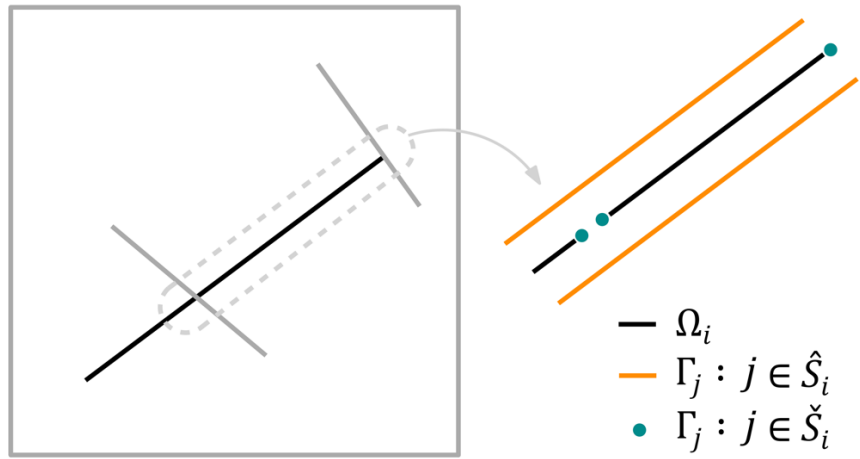

*Figure 3: Illustration of the interfaces neighboring a 1d fracture subdomain $\Omega_i$. Interfaces towards higher-dimensional neighbors are represented by the index set $\hat{S}_i$, while the corresponding index set for lower-dimensional neighbors is denoted $\check{S}_i$.*

## 3.1 Flow in fractured media

We first consider flow in mixed-dimensional geometries, where we assume a Darcy-type relation between flux and pressure gradients in all subdomains. The model has been presented several times before, see *e.g.* [2], [47]; our presentation therefore emphasizes implementational aspects within the framework presented in the previous sections. Moreover, of the model problems considered in this work, the flow problem by far has the best developed mathematical theory, and we therefore use this section to illustrate both variational and integral approaches to mixed-dimensional modeling and simulation.

### 3.1.1 Governing equations

To introduce the model, we first consider a domain with a single interface $\Gamma_j$, with neighboring subdomains $\Omega_h$ and $\Omega_l$, such that $\Omega_h$ is of one dimension higher than $\Gamma_j$ (which thus acts as an internal boundary) and $\Omega_l$ geometrically coincides with $\Gamma_j$ as illustrated in Figure 2. The flow model presented below has been studied *e.g.* in [39], [47]. We denote the flux on $\Gamma_j$ by $\lambda_j$, we can formally write $\lambda_j = \Pi_j^h tr\, q_h \cdot n_h$, with $n_h$ the unit normal of $\partial_j\Omega_h$ pointing from $\Omega_h$ to $\Omega_l$ and $tr$ a suitable trace operator from $\Omega_h$ to $\partial_j\Omega_h$. First consider the strong form of the Darcy problem stated for $\Omega_l$, which reads: find $(q_l, p_l)$ such that

$$q_l + \frac{\mathcal{K}_l}{\mu_l}\nabla p_l = 0,$$
$$\nabla \cdot q_l - \Xi_j^l \lambda_j = f_l, \tag{3.1}$$

where the differential operators are defined on the tangent space of $\Omega_l$ and $\Xi_j^l$ maps from $\Gamma_j$ to $\Omega_l$. We have indicated with $f_l$ a scalar source or sink term, $\mu_l$ is the fluid viscosity, while $\mathcal{K}_l$ represents the effective tangential permeability tensor, scaled by aperture, for more information see [48]. An analogous problem is written also for $(q_h, p_h)$, with the exception that $\Xi_j^h \lambda_j$ is mapped to a boundary condition on $\partial_j\Omega_h$,

$$q_h \cdot n_h|_{\partial_j\Omega_h} = \Xi_j^h \lambda_j. \tag{3.2}$$

The flux $\lambda_j$ is given by an interface condition of Robin-type on $\Gamma_j$ for $\Omega_h$, which reads

$$\lambda_j + \frac{\kappa_j}{\mu_j}\left(\Pi_j^l p_l - \Pi_j^h tr\, p_h\right) = 0,$$

(3.3)

where $\kappa_j$ indicates the normal effective permeability, and $\Pi_j^l$ and $\Pi_j^h$ the normal projection operators to $\Gamma_j$ from $\Omega_l$ and $\partial_j\Omega_h$, respectively. Equation (3.3) can be seen as a Darcy law in the normal direction associated to $\Gamma_j$. Several boundary conditions can be imposed on the external boundary of $\Omega_h$ and $\Omega_l$, for simplicity we limit ourselves to homogeneous pressure conditions (in the following). If $\Omega_l$ has a portion of the boundary which does not touch the external boundary, the so-called tip condition will be imposed being null flux. For notational convenience, we consider a unit viscosity for the remainder of this section; the viscosity is reintroduced in Section 3.2.

The extension to problems with many subdomains is now immediate: The interface still relates to its two neighboring subdomains, while for a subdomain $\Omega_i$ summation over all neighboring interfaces gives the problem: Find $(q_i, p_i)$ so that

$$q_i + \mathcal{K}_i \nabla p_i = 0,$$
$$\nabla \cdot q_i - \sum_{j \in \hat{S}} \Xi_j^i \lambda_j = f_i,$$
$$q_i \cdot n_i|_{\partial_j\Omega_i} = \Xi_j^i \lambda_j \quad \forall j \in \check{S}_i.$$

(3.4)In the case of $d = 0$, most of the above terms are void, and we are left only with the balance between the source term and fluxes from higher dimensions.

### 3.1.2 Subdomain variational and integral formulation

To move towards the numerical implementation, we introduce variational and integral formulations of the problem (3.4). Again, we focus on a domain with two subdomains and a single interface; the extension to several subdomains is straightforward.

Let the mortar variable be represented by $\lambda_j \in W(\Gamma_j) = L^2(\Gamma_j)$, in this setting we can define more precisely the projections $\Xi_j^l: L^2(\Gamma_j) \to L^2(\Omega_l)$ and $\Xi_j^h: L^2(\Gamma_j) \to L^2(\partial_j \Omega_h)$ that map interface fluxes to the neighboring subdomains. We note that the fluxes are mapped to the boundary of $\Omega_h$ but to the interior of $\Omega_l$, hence $\lambda_j$ acts as a boundary condition and a source term for the higher- and lower-dimensional subdomain, respectively. We also particularize projections from subdomains to the interface, defined as $\Pi_j^l: L^2(\Omega_l) \to L^2(\Gamma_j)$ and $\Pi_j^h: L^2(\partial_j \Omega_h) \to L^2(\Gamma_j)$; we shall comment on the implied $L^2$ regularity on $\partial_j \Omega_h$ below.

We first develop a mixed variational formulation of (3.4), introducing the following functional spaces

$$V(\Omega_h) = \{v \in H_{\nabla\cdot}(\Omega_h): tr\, v \cdot n_h \in L^2(\partial_j \Omega_h)\}, \quad V(\Omega_l) = H_{\nabla\cdot}(\Omega_l) \quad \text{and} \quad Q(\Omega_l) = L^2(\Omega_l).$$

Moreover, let the space $Q(\Omega_h)$ be a subspace of $L^2(\Omega_h)$ such that it is possible to define the operator $tr\, p_h$ with range at least in $L^2(\partial_j \Omega_h)$. It is well known that the trace cannot be defined for $L^2$-functions, however we note that, for example, the space $H^1(\Omega_h)$ fulfils the requirements. The extra request for $V(\Omega_h)$ on the interface is due to the low regularity of the trace on $H_{\nabla\cdot}(\Omega_h)$, which is related to the Robin-type nature of the coupling condition, see [2], [49].

The weak formulation of the mixed-dimensional Darcy problem reads: find $(q_h, p_h, q_l, p_l, \lambda_j) \in V(\Omega_h) \times Q(\Omega_h) \times V(\Omega_l) \times Q(\Omega_l) \times W(\Gamma_j)$ such that

$$\begin{aligned}
(\mathcal{K}_h^{-1} q_h, v)_{\Omega_h} - (p_h, \nabla \cdot v)_{\Omega_h} + (\Xi_j^h \lambda_j, tr\, v \cdot n_h)_{\partial_j \Omega_h} &= 0 & \forall v \in V(\Omega_h), \\
-(\nabla \cdot q_h, w)_{\Omega_h} &= -(f_h, w)_{\Omega_h} & \forall w \in Q(\Omega_h), \\
(K_l^{-1} q_l, v)_{\Omega_l} - (p_l, \nabla \cdot v)_{\Omega_l} &= 0 & \forall v \in V(\Omega_l), \\
-(w, \nabla \cdot q_l)_{\Omega_l} + (\Xi_j^l \lambda_j, w)_{\Omega_l} &= -(f_l, w)_{\Omega_l} & \forall w \in Q(\Omega_l), \\
(\kappa_j^{-1} \lambda_j, \mu)_{\Gamma_j} + (\Pi_j^l p_l, \mu)_{\Gamma_j} - (\Pi_j^h tr\, p_h, \mu)_{\Gamma_j} &= 0 & \forall \mu \in W(\Gamma_j).
\end{aligned}$$

(3.5)

Here $(\cdot,\cdot)_A$ is the $L^2$-scalar product on the set $A$. The problem is well posed as shown in [47]. We emphasize that, apart from the extra regularity assumptions on $\partial_j \Omega_h$, the variational formulation for the subdomains have the same structure as a fixed-dimensional problem.

Next, we state an integral formulation for the subdomain problems, expressed in primal form, that is, considering only the pressure variable. To that end, let $\omega_h \subset \Omega_h$ and $\omega_l \subset \Omega_l$ be subdomains (grid cells in the discrete setting) in the higher- and lower-dimensional domains. Likewise, $\omega_\gamma \subset \Gamma_j$ is a subdomain of the interface. In the view of subsequent considerations, we assume that the subdomains $\omega_h$ are non-overlapping and fully cover $\Omega_h$, and similarly for $\omega_l$ and $\omega_\gamma$. Additional requests on the shape regularity of $\omega_h$, $\omega_l$ and $\omega_\gamma$ depend on the numerical scheme. The integral formulation of (3.4) then reads: find $(p_h, p_l, \lambda_j)$ such that

$$\int_{\partial \omega_h \setminus \partial_j \Omega_h} \mathcal{K}_h \nabla p_h \cdot n_{\partial \omega_h} \, d\sigma + \int_{\partial \omega_h \cap \partial_j \Omega_h} \Xi_j^h \lambda_j \, d\sigma = \int_{\omega_h} f_h \, dx \qquad \omega_h \subset \Omega_h,$$

$$\int_{\partial \omega_l} \mathcal{K}_l \nabla p_l \cdot n_{\partial \omega_l} \, d\sigma - \int_{\omega_l} \Xi_j^l \lambda_j \, dx = \int_{\omega_l} f_l \, dx \qquad \omega_l \subset \Omega_l,$$

$$\int_{\omega_\gamma} \kappa_j^{-1} \lambda_j \, dx + \int_{\omega_\gamma} \Pi_j^l p_l \, dx - \int_{\omega_\gamma} \Pi_j^h tr \, p_h \, dx = 0 \qquad \omega_\gamma \subset \Gamma_j,$$

(3.6)

where we have indicated with $dx$ and $d\sigma$ the infinitesimal measure for equi-dimensional and one co-dimensional integrals, respectively, with respect to the considered cell dimension. The vector $n_{\partial \omega}$ denotes the outward unit normal of $\omega$. The equations are written on all subdomains $\omega$ and the global problem is given once the continuity of normal fluxes is imposed on each $\partial \omega$.

We make two remarks related to the discretization of the above equations. First, both the variational and integral formulations are very close to the corresponding fixed-dimensional problems, thus there is considerable scope for reuse of existing software. Second, seen from the interface, the subdomain discretization acts as an unspecified Neumann-to-Dirichlet map that converts the interface fluxes into pressures to be projected to the interface. The formulation is independent of the actual discretization on the subdomains, and there is no requirement the same discretization be used on the two neighboring subdomains (of an interface). For more information on the formulation, confer [39].

### 3.1.3 Implementation

From the variational and integral formulations stated above, we see that for a discretization on a generic subdomain $\Omega_i$ to interact with the interface problem, we need to provide operators which:

1) Handle Neumann boundary data on the form $\Xi_j^i \lambda_j$, for all interfaces $\Gamma_j$ where $\Omega_i$ is the higher-dimensional neighbor.
2) Handle source terms $\Xi_j^i \lambda_j$ from interfaces $\Gamma_j$ where $\Omega_i$ is the lower-dimensional neighbor.
3) Provide a discrete operator $tr \, p_i$ so that $\Pi_j^i$ can project the pressure trace from $\partial_j \Omega_i$ to $\Gamma_j$ where $\Omega_i$ is the higher-dimensional neighbor.
4) Provide a pressure $p_i$ so that $\Pi_j^i$ can project the pressure to all $\Gamma_j$ where $\Omega_i$ is the lower-dimensional neighbor.

Of these, all but the third operation is readily available in any reasonable implementation of a discretization scheme for elliptic equations. For the discrete pressure trace there is some room for interpretation; the simplest approach is to associate the trace with the pressure in cells immediately next to the interface. Higher order trace reconstruction operators, utilizing the construction of the discretization scheme at hand, are also possible; our implementation of finite volume methods for flow allows for sub-cell variations in pressure so that the discrete pressure at $\partial_j \Omega_i$ differs from the cell center pressure closest to the boundary.

It is instructive to write out the structure of the coupled system for our case with two subdomains $\Omega_h$ and $\Omega_l$ separated by an interface $\Gamma_j$. Denote by $y_h$, $y_l$ and $\xi_j$ the vectors of discrete unknowns in $\Omega_h$, $\Omega_l$ and on $\Gamma_j$, respectively. As we make no assumptions that the same discretization scheme is applied in both subdomains, these may contain different sets of unknowns. Specifically, the unknown can be cell center pressures only, or cell center pressure and face fluxes, depending on the discretization scheme applied. The discrete coupled system can then be represented on the generic form

$$\begin{pmatrix} A_h & 0 & N_h \Xi_j^h \\ 0 & A_l & S_l \Xi_j^l \\ -\Pi_j^h P_h & \Pi_j^l P_l & D_j \end{pmatrix} \begin{pmatrix} y_h \\ y_l \\ \xi_j \end{pmatrix} = \begin{pmatrix} f_h \\ f_l \\ 0 \end{pmatrix}. \tag{3.7}$$

Here, $A_h$ and $A_l$ are the fixed-dimensional discretizations on the subdomains, $N_h$ is the discretization of Neumann boundary conditions on $\Omega_h$, and $S_l$ is the discretization of source terms in $\Omega_l$. Furthermore, $P_h$ provides a discrete representation of the pressure trace operator on $\Omega_h$ and $P_l$ gives the pressure unknowns in $\Omega_l$; the latter is an identity operator for the integral formulations presented on primal form and strips away flux unknowns in the dual formulation. Finally, $D_j$ is the discretization of (3.3). In accordance with the second constraint on mixed-dimensional modeling discussed in Section 2.2, there is no direct coupling between $\Omega_h$ and $\Omega_l$. Global boundary conditions are left out of the system; as a technical detail we note that for some discretization schemes, *e.g.* multi-point flux approximation (Mpfa) methods, the global boundary conditions can also give a contribution to the right-hand-side of the interface equation.

The form (3.7) suggests an implementation strategy, based on the graph representation of the mixed-dimensional domain, which also exploits reuse of software for fixed-dimensional problems: On the graph nodes, that is the subdomains, the pressure equation is discretized as if it were a fixed-dimensional problem. The interface law is discretized by traversal of the graph edges; this operation will communicate with the discretizations in the neighboring subdomains to obtain the terms represented in the last column and row of (3.7).

The PorePy implementation of the above method represents the mortar variable by piecewise constant functions. Due to the decoupling, there is no requirement that the same numerical method be used on all subdomains, and indeed PorePy gives complete flexibility in this respect by an implementation of the coupling structure (3.7) which is independent of the individual subdomain discretizations. PorePy offers four discretization schemes for the flow problem: Lowest order Raviart-Thomas mixed finite elements combined with a piecewise constant pressure approximation (RT0-P0) [49], the lowest order mixed virtual element method (Mvem) [50], [51], and two finite volume schemes: the two- and multi-point flux approximations (Tpfa, Mpfa) [52]–[54]. Our implementation for the coupled mixed-dimensional problem relies on the analysis carried out in [39], which provides a theoretical background to obtain a stable global scheme.

### 3.2 Fully coupled flow and transport

We next turn to simulation of fully coupled flow and transport, as an example of a multiphysics problem with variable couplings within and between subdomains. We consider the mixing of two incompressible and miscible species of different viscosities. We put emphasis on the modeling of the mixed-dimensional dynamics and discuss some implementation aspects. For discretization, we limit ourselves to finite volume methods for the problem written in primal form. We present the governing equations on integral form only with details on how to handle the advective part at the interfaces.

### 3.2.1 Continuous formulation

We start by considering a single subdomain $\Omega_i$. Denote the pressure in a subdomain $\omega_i \subset \Omega_i$ by $p_i$. We represent the species evolution by the mass concentration $c_i$ in $\omega_i$. By the incompressibility of the fluids, the conservation of total mass within $\omega_i$ can be written as

$$\int_{\partial \omega_i} q_i \cdot n_{\partial \omega_i} \, d\sigma = \int_{\omega_i} f_i \, dx, \tag{3.8}$$

where $f_i$ represents the total volumetric sources and sinks. The Darcy flux $q_i$, depends on both pressure and mass concentration, via the fluid viscosity $\mu = \mu(c_i)$, and is given by

$$q_i + \frac{\mathcal{K}_i}{\mu(c_i)} \nabla p_i = 0. \tag{3.9}$$

Here $\mathcal{K}_i$ denotes the effective tangential permeability of $\Omega_i$. Conservation of mass for each species is expressed by the equation

$$\int_{\omega_i} \phi_i \frac{\partial c_i}{\partial t} \, dx + \int_{\partial \omega_i} w_i \cdot n_{\partial \omega_i} \, d\sigma = \int_{\omega_i} g_i \, dx \tag{3.10}$$

Here, $\phi_i$ represents the effective porosity, $g_i$ denotes sources and sinks for the species, and the flux $w_i$ is composed by a diffusive and an advective term

$$w_i + \mathcal{D}_i \nabla c_i - c_i q_i = 0, \tag{3.11}$$

where $\mathcal{D}_i$ is the effective diffusivity of $\Omega_i$. We note that the equations are coupled via the concentration dependency of viscosity and the presence of the Darcy flux in the advective transport.

The interaction between two neighboring subdomains $\Omega_h$ and $\Omega_l$ again goes via the common interface $\Gamma_j$. The total flux over $\Gamma_j$, denoted by $\lambda_j$, is given by (3.3), where the interface viscosity $\mu_j$ is modeled as a function of the mean of the concentrations on the two sides,

$$\mu_j = \mu_j \left( \frac{\Pi_j^l c_l + \Pi_j^h tr \, c_h}{2} \right). \tag{3.12}$$

Mass flux over $\Gamma_j$ is again governed by an advection-diffusion relation: The diffusion term $\beta_j$ is, in analogy with the corresponding term for the Darcy flux, given by

$$\beta_j + \delta_j \left( \Pi_j^l c_l - \Pi_j^h tr \, c_h \right) = 0, \tag{3.13}$$

with $\delta_j$ representing the effective diffusivity over the interface $\Gamma_j$. For the advective term $\eta_j$, we introduce an upstream-like operator based on the Darcy interface flux:

$$Up(c_h, c_l; \lambda_j) = \begin{cases} \Pi_j^h tr\, c_h, & \text{if } \lambda_j \geq 0 \\ \Pi_j^l c_l, & \text{if } \lambda_j < 0. \end{cases} \tag{3.14}$$

With this, the advective interface flux $\eta_j$ is given by the relation

$$\eta_j - \lambda_j Up(c_h, c_l; \lambda_j) = 0. \tag{3.15}$$

What remains in the problem formulation is to introduce the coupling terms in the subdomain equations and introduce global boundary conditions. As all the interface variables are fluxes, their treatment is analogous to that discussed for the flow problem in Section 3.1. For $\Omega_h$, the interface fluxes enter as flux boundary conditions for the total mass flux ($\lambda_j$) and mass concentration ($\beta_j, \eta_j$) conservation equations, while for $\Omega_l$, the fluxes enter as corresponding source terms. Finally, global boundary conditions are imposed in the standard way for elliptic and advection-diffusion problems, see *e.g.* [55]. With few modifications, our formulation can handle a purely advective problem, like transport of a passive scalar. In this case the elliptic operators (interface and mortar law included) are not considered.

The equations (3.8)-(3.15) define the governing equations in all subdomains and on all interfaces. The only exception is 0d domains, wherein the fluid mass and concentration fluxes are void, and the governing equations simply balance the fluxes of neighboring interfaces with possible source terms in the point domain.

### 3.2.2 Implementation

The equations are discretized with finite volume methods, by letting $\omega_i$ represent a computational cell. To discretize the flux expression, we apply single point upstreaming for the advective flux [56] and Mpfa for the diffusive terms both in (3.8)-(3.9) and (3.10)-(3.11), as described in Section 3.1. These operations can be carried out independently on individual subdomains, and they can readily reuse existing implementations for fixed-dimensional problems.

To discretize the interface laws, we need projection operators for scalar quantities (pressure and mass concentration) from subdomains to interfaces, and projections of fluxes from interfaces back to subdomains. For the diffusive fluxes, the treatment is identical to that described in Section 3.1 for the elliptic equation. Similarly, the advective terms will appear respectively as flux boundaries and source terms for the higher- and lower-dimensional neighbors of an interface. Handling of these extra terms should be straightforward in any existing code for fixed-dimensional problems.

Finally, we note that governing equations are non-linearly coupled via the viscosity and the presence of the Darcy flux in the advective transport terms. A non-linear solver is therefore needed. Within PorePy, this is most easily handled by an automatic differentiation module; for an example of how this is applied, confer the supplementary material for the simulation script used in Section 5.1.

### 3.3 Poromechanical fracture deformation by contact mechanics

Our final set of model equations considers poroelastic deformation of a fractured medium, where the fractures may open and / or undergo abrupt slip if the frictional forces are insufficient to withstand tangential tractions on the fracture surface. This process is important in applications such as geothermal energy extraction and $CO_2$ storage. Moreover, modeling of the process leading up to and under sliding is non-trivial due to i) the coupled poreelastic processes, ii) heterogeneous modeling equations between subdomains, iii) the need to use non-standard constitutive laws to relate primary variables during sliding, iv) non-smooth behavior of the constitutive laws in the transition from sticking to sliding of a fracture. Modeling of this process is an active research field, see *e.g.* [57]–[59], and so this is an example where the availability of a flexible prototyping framework for this research is extremely useful.

Herein, we present a set of governing equations which borrows modeling concepts from contact mechanics to describe the sliding problem. Our formulation has no notion of a displacement inside the fracture, instead the fracture deformation is described by the displacement jump over the fracture surface. We show how the model is naturally formulated and implemented in our mixed-dimensional framework, by defining displacement variables on the matrix-fracture interface.

### 3.3.1 Governing equations

As modeling of deformation of intersecting fractures is non-trivial, we limit our exposition to media with non-intersecting fractures. Flow and deformation in the rock matrix, represented by the subdomain $\Omega_h \in \mathbb{R}^n$, are then governed by Biot's equations for poroelasticity [60]

$$\nabla \cdot (\mathcal{C}_h \nabla_s u_h - \alpha_h p_h I) = b_h, \tag{3.16}$$

$$\alpha_h \frac{\partial (\nabla \cdot \dot{u}_h)}{\partial t} + \theta_h \frac{\partial p_h}{\partial t} - \nabla \cdot \left(\frac{\mathcal{K}_h}{\mu_h} \nabla p_h\right) = f_h.$$

Here, the first equation represents conservation of momentum, with the acceleration term neglected, while the second equation expresses conservation of mass. The primary variables are the displacements, $u_h$, and the fluid pressure $p_h$. The stiffness matrix $\mathcal{C}_h$ can for linear isotropic media be expressed purely in terms of the first and second Lamé parameters, and the stress can be computed as $\sigma_h = \mathcal{C}_h \nabla_s u_h$, where $\nabla_s$ is the symmetric gradient. Furthermore, $\alpha_h$ is the Biot constant, $I$ the second order identity tensor, $b_h$ denotes body forces, $\theta_h$ the effective storage term, $\mathcal{K}_h$ the permeability and $\mu_h$ the viscosity. We also assume boundary conditions are given on the global boundary.

Next, consider an interface $\Gamma_j$ between the higher-dimensional subdomain $\Omega_h$ and the lower-dimensional domain $\Omega_l$. Denote the displacement variable on $\Gamma_j$ by $u_j$. We emphasize that $u_j$ is a vector in $\mathbb{R}^n$, that is, it represents the displacement in both the tangential and normal direction of $\Omega_l$. We will require continuity between $u_h$ and $u_j$, expressed as $\Pi_j^h tr\, u_h = u_j$, where we recall that the trace operator maps to $\partial_j \Omega_h$, the part of the boundary of $\Omega_h$ that coincides with $\Gamma_j$. We also need to introduce the jump in displacement, $[\![u_j]\!]$, between the two interfaces on opposing sides of $\Omega_l$, see Figure 4 for an illustration. The jump is decomposed into the tangential jump $[\![u_j]\!]_\tau$ and the normal jump $[\![u_j]\!]_n$.

The mechanical state in $\Omega_l$ is described by the contact pressure $\sigma_l$, which again is a vector in $\mathbb{R}^n$, with tangential and normal components $\sigma_{l,n}$ and $\sigma_{l,t}$, respectively. Our model also includes fluid flow in the fracture $\Omega_l$, which is governed by conservation of mass

$$\frac{\partial}{\partial t}\left(a([\![u_j]\!])\right) + \theta_l \frac{\partial p_l}{\partial t} - \nabla \cdot (\mathcal{K}_l \nabla p_l) - \Xi_j^l \lambda_j = f_l.$$

(3.17)

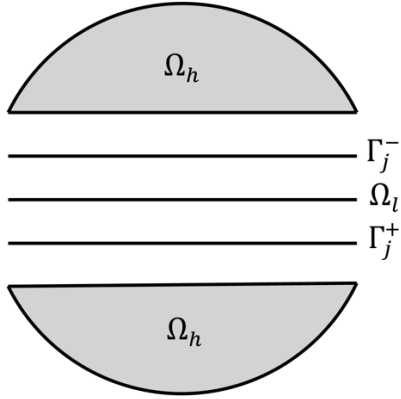

Figure 4: Illustration of a lower-dimensional domain, $\Omega_l$, that has two interfaces, $\Gamma_j^-$ and $\Gamma_j^+$, with a higher dimensional domain, $\Omega_h$. The inerfaces on opposing sides are identified by a positive and negative index.

Here, the time derivative of the aperture $a(\llbracket u_j \rrbracket)$ represents changes in the available volume due to changes in the displacement jump; in general this can be both reversible changes due to normal displacement of the fracture, and irreversible changes caused by shear dilation [3]. In the following, we only consider the normal part, *i.e.* $a(\llbracket u_j \rrbracket) = \llbracket u_j \rrbracket_n$. Similarly, it can be of interest to consider changes in the permeability $K_l$ due to changes in aperture, although we shall not consider this topic herein. Finally, $\theta_l$ denotes the effective storage term. As in the previous sections, the relation between the fluid pressures in $\Omega_h$ and $\Omega_l$ is governed by a flux law of the type (3.3), and we let $\lambda_j$ denote the interface variable that represents fluid flux between the domains.

The relation between $\sigma_l$ and $\llbracket u_j \rrbracket$ is modeled by borrowing techniques from contact mechanics as summarized here; for a full discussion see [61]. Balance of forces between the poroelastic stress in $\Omega_h$ and the contact pressure in $\Omega_l$ is expressed as

$$\Pi_j^h\, n_h \cdot (\sigma_h - \alpha_h p_h I) = \Pi_j^l \sigma_l - (\Pi_j^h\, n_h)(\Pi_j^l \alpha_l p_l) \quad \text{on } \Gamma_j. \tag{3.18}$$

In the direction normal to $\Omega_l$, the contact stress is zero only when the displacement jump is nonzero, that is

$$\llbracket u_j \rrbracket_n \leq 0, \quad \sigma_{l,n} \leq 0, \quad \llbracket u_j \rrbracket_n \sigma_{l,n} = 0. \tag{3.19}$$

The motion in the tangential direction is controlled by the ratio between the tangential force $\sigma_{l,\tau}$ and the maximum available frictional force $F\sigma_{l,n}$, where $F$ is the friction coefficient. The time derivative of the displacement jump is zero until the frictional force is overcome; for larger tangential forces, the derivative of the displacement jump and tangential force are parallel:

$$\begin{cases} ||\sigma_{l,\tau}|| \leq -F\sigma_{l,n}, \\ ||\sigma_{l,\tau}|| < -F\sigma_{l,n} \rightarrow \llbracket \dot u_j \rrbracket_\tau = 0, \\ ||\sigma_{l,\tau}|| = -F\sigma_{l,n} \rightarrow \exists \alpha \in \mathbb{R},\ \sigma_{l,\tau} = -\alpha^2 \llbracket \dot u_j \rrbracket_\tau, \end{cases} \tag{3.20}$$

where $||\cdot||$ represents the Euclidean norm and $\llbracket \dot u_j \rrbracket_\tau$ represents the sliding velocity. We emphasize that the contact conditions are formulated in terms of the contact stress $\sigma_l$, with no contribution from the pressure $p_l$.

### 3.3.2 Implementation in mixed-dimensional framework

It is instructive to discuss implementation of poroelastic contact mechanics within our mixed-dimensional modeling concept, starting from an existing implementation of poroelasticity in the matrix domain. This is a relevant case for many research codes, in particular the PorePy implementation for this problem was extended from a finite volume method, the multipoint stress approximation (Mpsa), originally developed for elastic and poroelastic deformation [62]–[64]. Below, we follow the equations presented above and identify variables and equations to be introduced.

First, the variables $p_l$ and $\lambda_j$ representing respectively fluid pressure in the fracture $\Omega_l$ and the fluid flux on the interface $\Gamma_j$ between $\Omega_h$ and $\Omega_l$, are introduced as discussed in detail in Section 3.1. This implies that $\partial_j \Omega_h$ is a Neumann boundary for fluid flow in $\Omega_h$. We reiterate that the couplings introduced by this approach is standard for any discretization scheme for single-phase flow.

Second, the displacement in $\Omega_h$ must be coupled to the mortar displacement $u_j$ on $\Gamma_j$. This is achieved by letting $\partial_j \Omega_h$ be a Dirichlet boundary, so that the condition $\Pi_j^h tr\, u_h = u_j$ can be enforced by the imposition of a boundary condition. The poroelastic stress at $\partial_j \Omega_h$ is computed from variables in $\Omega_h$ and on $\Gamma_j$, according to the discretization scheme applied in $\Omega_h$. We note that mapping of displacements $u_j$ onto the boundary $\partial_j \Omega_h$, and later stresses from $\partial_j \Omega_h$ to $\Gamma_j$ requires vectorized versions of the projection operators discussed in Section 2.4; this is a straightforward extension.

Finally, the implementation must discretize the stress continuity as expressed by (3.18), and the relation between displacement jumps $[\![u_j]\!]$ and contact pressure $\sigma_l$ (3.19)-(3.20). Stress continuity is enforced by projecting the discrete representation of the poroelastic stress on $\partial_j \Omega_h$ onto $\Gamma_j$, similarly projecting the discrete quantity $\sigma_l - p_l n$ from $\Omega_l$ to $\Gamma_j$ and enforcing equality. The contact conditions are discretized by projecting $[\![u_j]\!]$ onto $\Omega_l$, and then discretizing Equations (3.19) and (3.20). This is a non-linear term, in that the relation between $[\![u_j]\!]$ and $\sigma_j$ depends on whether the fracture is open, sticking or slipping. In our implementation we use a semi-smooth Newton method to deal with the discontinuities in the solution, for details we refer to [61], [65].

As a final remark on data structures, we note that the full discrete system is rather complex, with different governing equations in different subdomains, and non-trivial couplings between variables that live on different grids. As illustrated by the run scripts in the supplementary material, the mixed-dimensional grid structure and modeling concept break the implementation into manageable parts. Moreover, due to the strong modularization of the model and implementation, experimentation with model variations etc. is handled with minimal needs for adjustments.

## 4 Validation

To validate our modeling framework and its implementation in PorePy, we consider three test cases: A benchmark for flow problems in 2d fractured media, Mandel's problem for poroelasticity, and Sneddon's problem for fracture deformation in elastic media. Together, these cases probe a wide range of the capabilities of the modeling framework and its PorePy implementation, including discretization schemes, multiphysics problems and time-dependent problems. The cases thus supplement previous testing of PorePy, reported in [39], [66]–[68]. The supplementary material provides detailed setups, including parameters for all simulations in Section 4 and 5. Scripts that reproduce all results reported herein can be accessed at [69]; see that reference or the supplementary material for install instructions.

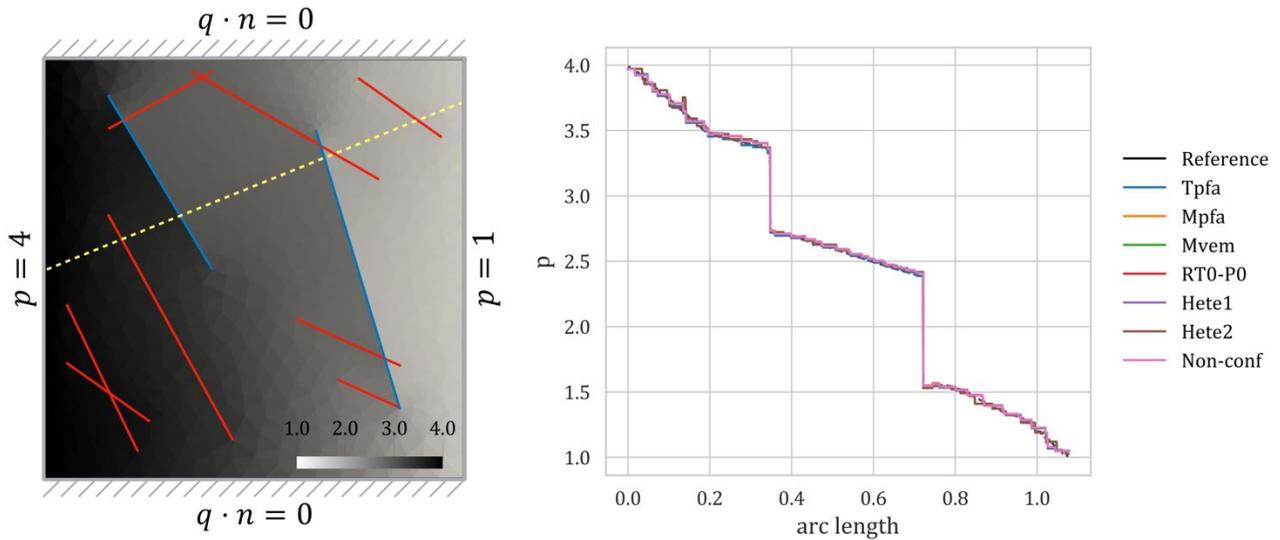

*Figure 5: Left: A solution obtained with Mpfa on the coarsest mesh, showing the fracture network and the problem setup. The red lines represent conductive fractures whereas the blue lines are blocking fractures. The yellow line indicates the line of the pressure profile. Right: Pressure profiles for the discretization schemes used in the validation.*

## 4.1 Flow in 2d fractured porous media

To validate the mixed-dimensional flow discretization, we consider Benchmark 3 of [70], which describes an incompressible single-phase flow problem in a fractured domain. The fracture network contains intersecting and isolated fractures; see Figure 5 for an illustration of the domain together with the pressure solution for the Mpfa discretization. The network contains both highly conductive and blocking fractures, see the supplementary material for details. The normal permeability in the fracture intersections is given by the harmonic average of the permeabilities of the intersecting fractures as suggested in [70], [71].

The aim of this case is twofold - we benchmark our code to well-established methods in the literature and present the full capability that our abstract structure can handle. For the latter, we consider four groups of discretization schemes and simulation grids: First, three homogeneous (the same for all the subdomains) discretizations: Tpfa, Mpfa and RT0-P0. Second, a case with the Mvem, where the cells of the rock matrix are constructed by a clustering procedure starting from a more refined simplicial grid, see [67] for details. Third, two heterogeneous discretizations where RT0-P0 and Mvem for the rock matrix are combined with Tpfa for the fractures (labeled Hete1 and Hete2, respectively). Fourth, a case where the fracture grid is twice as fine as the matrix grid, with the mortar grids non-conforming to the surrounding grids (Non-conf). In this case we consider the RT0-P0 scheme. We use simplex grids in all cases that do not involve Mvem.

Figure 5 shows the domain with fractures, boundary conditions and a representative numerical solution. The figure also depicts a plot of the pressure along the line $(0,0.5) - (1,0.9)$, with the reference solution (equi-dimensional problem computed on a much finer grid) colored in black. We observe good agreement between the solutions obtained in PorePy and the reference solution. We also consider a sequence of grids to compute the error relative to the reference solution, as done in the original benchmark. Figure 6 shows the decay of the normalized $L^2$ error for the rock matrix and the union of the fracture subdomains. In the former, we notice a first order of convergence for all the considered methods. The convergence rate for the fracture subdomains is sublinear, as was also observed in the original benchmark, see [70].

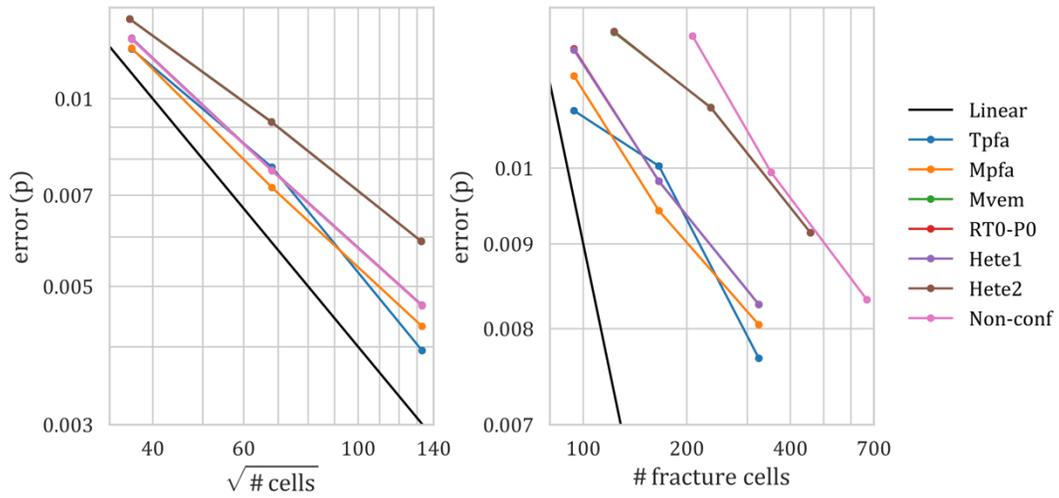

Figure 6: Left: Convergence of the pressure unknown for the matrix subdomain. Right: Convergence for the pressure unknown for the fracture subdomains.

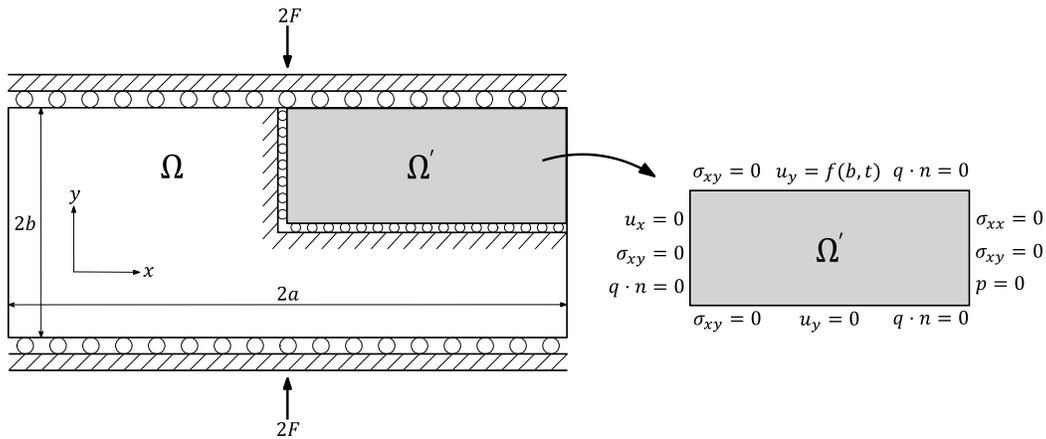

Figure 7: Mandel's problem. Left: Schematic representation of the full and positive quarter domains, $\Omega$ and $\Omega'$. Right: Quarter domain showing the boundary conditions.

## 4.2 Mandel's problem in poroelasticity

The next test case considers a poroelastic material, with a setup defined by Mandel's problem [72], [73], for which an analytical solution is available. While the problem geometry does not include lower-dimensional objects, the case tests the implementation of the poroelastic code and shows the framework's flexibility to deal with coupled problems and time-dependent mixed boundary conditions. The original problem consists of an isotropic poroelastic slab of width $2a$ and height $2b$ sandwiched by two rigid plates (Figure 7). Initially, two compressive constant loads of intensity $2F$ are applied to the slab at $y = \pm b$. At $x = \pm a$, fluid is free to drain, and edges are stress free. Gravity contributions are neglected.

This problem is modeled using the quasi-static Biot equations, as presented in Section 3.3.1. Exploiting the symmetry of the problem, we focus on the positive quarter domain $\Omega'$, rather than the full domain $\Omega$, see Figure 7 for an illustration, and for boundary conditions. Note that the vertical displacement at the top of the domain is time-dependent and given by the exact solution, see [74].

The simulation parameters were taken from [75], see also the supplementary material for details. The coupled problem is discretized in space using Mpsa/Mpfa for the mechanics and flow, respectively. For the time discretization we use implicit Euler. The computational mesh is unstructured and composed by 622 triangular elements. The results are shown in Figure 8 in terms of dimensionless quantities and are in good agreement with [75] for both pressure and displacement.

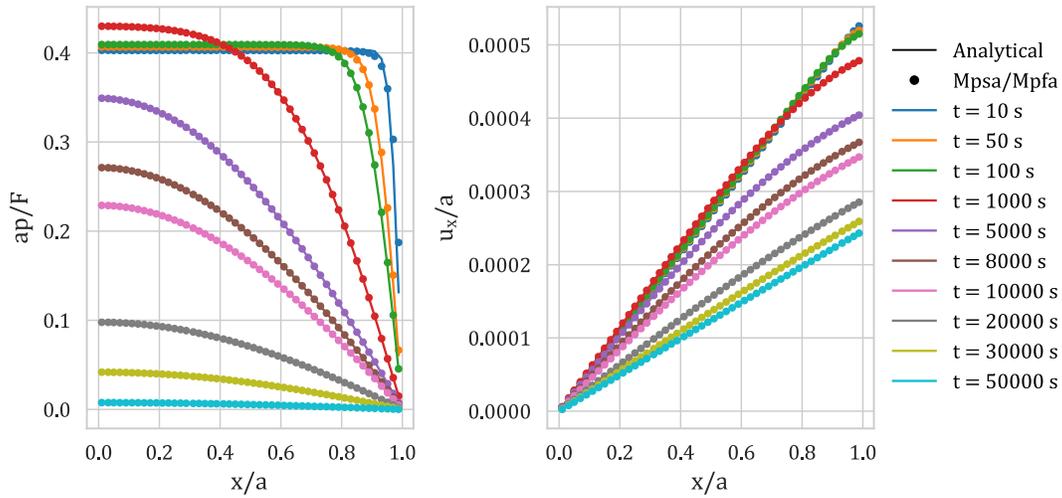

Figure 8: Solutions to Mandel's problem: Dimensionless pressure (left) and horizontal displacement (right) profiles for several times.

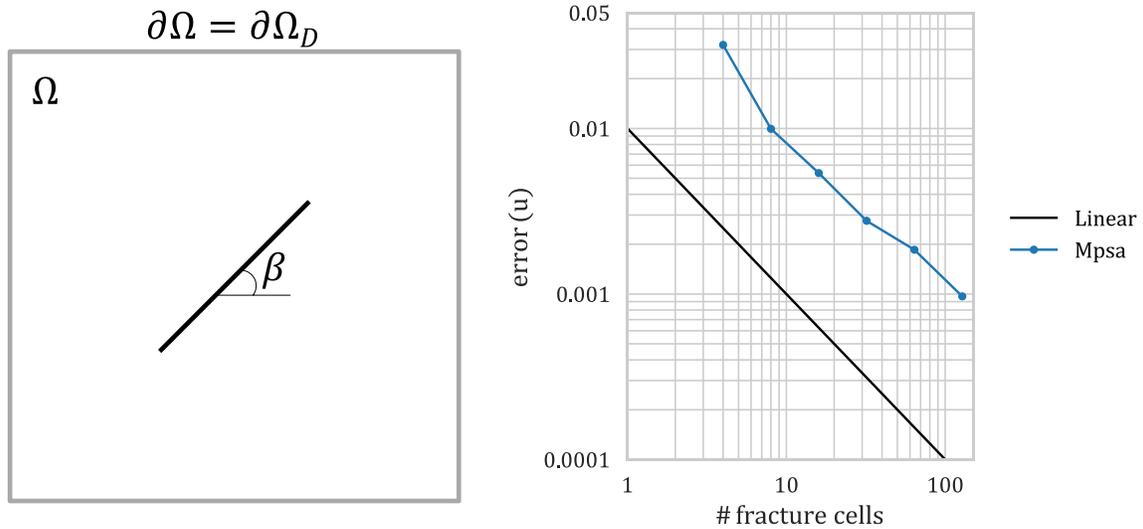

Figure 9: Setup and convergence of Sneddon's problem Left: Schematic representation of the domain. Right: Average convergence behavior of the relative normal displacement along the fracture. Each dot corresponds to the average of 140 simulations.

## 4.3 Sneddon's problem of fracture deformation

In this example, a square domain with a single fracture located in the middle is considered. The fracture forms an angle β with the horizontal direction (see Figure 9) and is subjected to a constant pressure $p_0$ acting on its interior. This pressure can be interpreted as a pair of normal forces acting on either side of the fracture. An analytical solution for the relative normal displacement along the fracture was derived by Sneddon [76] for an infinite domain, and has the following form:

$$[\![u_j]\!]_n(d_f) = \frac{(1-\nu)p_0 L}{G}\sqrt{1-\frac{d_f^2}{\left(\frac{L}{2}\right)^2}}$$

(4.1)

where $\nu$ and $G$ are the Poisson's ratio and shear modulus, respectively, $L$ is the fracture length, and $d_f$ denotes the distance from the center of the fracture.

In our calculations, the conditions of infinite domain are replaced with a Dirichlet boundary, where the prescribed displacement is set equal to the analytical solution calculated using the procedure illustrated in [77]. The accuracy of the numerical solution is very sensitive to the discretization, specifically the cell configuration, at the fracture tips [58]. To reduce the dependency on specific grid realizations, the values of the numerical solution reported in Figure 9 are the average of a group of 20 x 7 = 140 computations per level of grid resolution, with 7 different fracture angles $\beta$ in the range $0° - 30°$ and 20 grid realizations per fracture. With six levels of grid refinement, the full study contains 20 x 7 x 6 = 840 simulations. Figure 9 summarizes the results in the form of the error in relative normal displacement between the analytical solution (4.1) and the numerical solution as a function of the fracture resolution, i.e. number of fracture elements. The method provides first-order convergence on average.

## 5 Applications: Multiphysics simulations

Having established the accuracy of PorePy for central test cases that involve mixed-dimensional geometries, we go on to present two multiphysics cases of high application relevance: A fully coupled flow and transport problem, and fracture reactivation caused by fluid injection. The motivation for the simulations is to illustrate further capabilities of the modeling framework and its PorePy implementation, including simulations on complex 3d fracture networks, automatic differentiation applied to non-linear problems, non-matching grids, and simulation of fracture deformation in a poromechanical setting.

### 5.1 Fully coupled flow and transport

This example has two main purposes. First, we consider a non-linear coupled flow and transport process as described in Section 3.2. We apply the automatic differentiation functionality in PorePy to obtain the Jacobian of the global system of equations, which is then used in a standard Newton method to solve the non-linear problem. Second, we illustrate the flexibility of the mixed-dimensional approach by using non-matching meshes on a relatively complex fracture network.

We consider the injection of a highly viscous fluid into a domain initially filled with a less viscous fluid. The two fluids are miscible, with their distribution described by the mass concentration $c \in [0, 1]$, and with a viscosity ratio of the two fluids given by $\mu(c) = \exp(c)$. In the parameter regime studied in this example, the transport in the fractures is advection dominated, while the transport in the rock matrix is dominated by diffusion; see the supplementary material for the details about the parameters. We remark that PorePy has also been applied to study unstable displacement in 2d domains, see [78] for details.

The mixed-dimensional domain considered in this example consists of one 3d domain, 15 2d fracture domains, 62 1d domains and 9 0d domains. On this geometry, two computational grids are constructed: The first has matching grids in all dimensions, with in total 20812 cells, out of which 16766 are 3d cells and 3850 are 2d fracture cells. The second mixed-dimensional grid has a 3d grid identical to the first grid, whereas the lower-dimensional objects are assigned refined grids with in total 13839 2d fracture cells, thus the 3d-2d interfaces have non-matching grids. The combination of the non-linearity and the non-matching grids provides a challenging test for the robustness of the PorePy implementation of subdomain couplings and provides an illustration of the framework's flexibility.

Figure 10 shows the average concentration profile in the fractures for the two meshes. There are no significant differences in the average concentration profiles in the two cases, indicating the stability of the implementation of the non-matching case. Figure 11 shows a snapshot of the concentration in the fractures and the rock matrix at time $t = 20$. The diffusive front in the rock matrix has only moved a few grid cells at the break-through, however, due to the diffusion and advection from the fractures to the rock matrix, the concentration has increased in considerable parts of the rock matrix. We observe no irregularities for the solution produced on the non-matching grid, confirming PorePy's ability to deal with non-standard grid couplings also for challenging physical regimes.

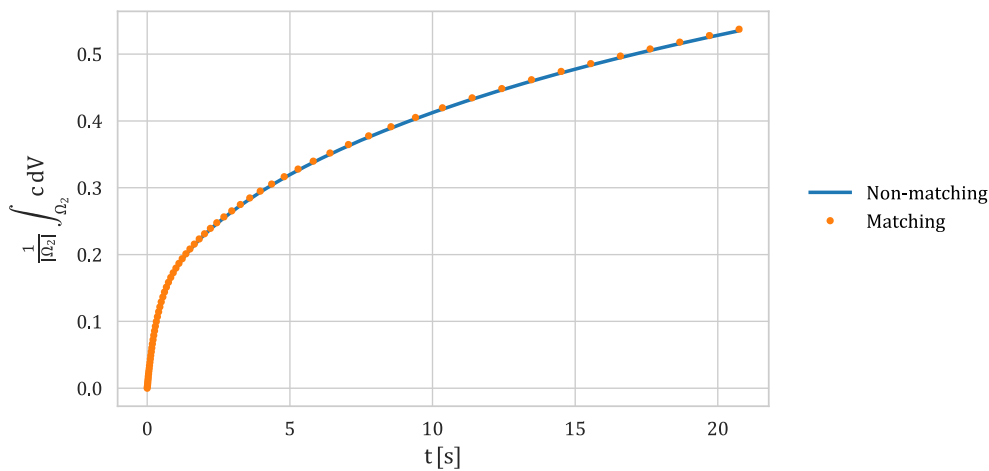

*Figure 10. Fully coupled flow and transport: Comparison of average concentration in the fracture network for a simulation with matching meshes and a simulation with non-matching meshes.*

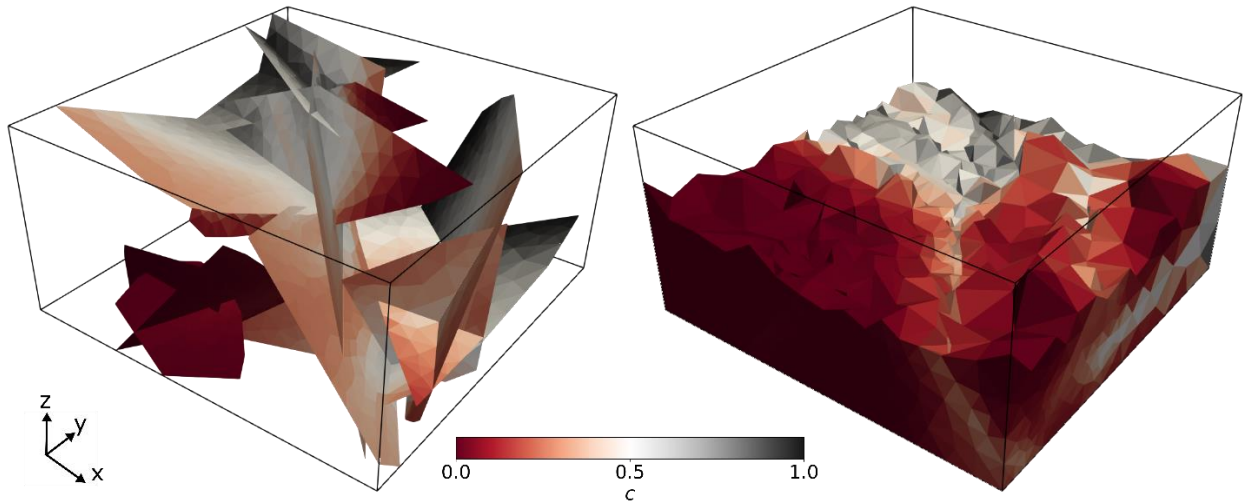

*Figure 11. Fully coupled flow and transport: Concentration in the fractures (left) and in the rock matrix (right) for the coupled flow and transport problem given in Section 3.2 at the end time of the simulation (t=20). In the right figure the rock matrix domain is cropped, and the fractures removed to reveal the concentration inside the domain. The black lines indicate the domain boundary. A non-matching mesh is used where the mesh in the fractures are much finer than the mesh in the rock matrix. The increase in the concentration in the rock matrix is mainly due to advection and diffusion from the fractures.*

## 5.2 Poroelasticity and fracture deformation

The final example aims at demonstrating the modelling framework's and PorePy's applicability to non-standard combinations of physical processes in different domains and thereby its potential for method development and prototyping. With the critical events taking place on individual fractures as a result of processes in the rock matrix, it also serves as an example of the importance of incorporating dynamics of both the matrix and explicitly represented fractures, as done in DFM models.

We consider a reservoir of idealized geometry containing three fractures numbered from 1 through 3, whereof the first contains an injection well, see Figure 12. On this geometry, we solve the governing equations presented in Section 3.3. We impose injection over a 25-day period and an anisotropic background stress regime, producing a scenario well suited to demonstrate different fracture dynamics. We investigate the dynamics both during a 25-day injection phase, and during the subsequent 25-day relaxation phase, at the end of which the pressure has almost reached equilibrium once more. The full set of parameters may be found in the supplementary material.

The dynamics on the fractures throughout the simulation are summarized in Figure 12, while the spatial distribution of the fracture displacement jumps at the end of the injection phase is shown in Figure 13. During the injection phase, there are tangential displacement jumps on all three fractures, appearing first on the favorably oriented fractures 2 and 3, and then on fracture 1 (injection). Normal displacement jumps appear on fracture 1 along with the tangential jumps, and on fracture 3 somewhat later. On fracture 2, which is located furthest away from the injection point, no normal displacement jumps appear. During the relaxation phase, we note that while the normal displacement jumps vanish, all tangential jumps remain due to the friction. We also observe a slight increase in tangential displacement jumps on fractures 2 and 3 at the time of shut-in, as the normal jump on fracture 1 vanishes.

The example demonstrates how modeling of complex coupled processes in great detail is possible through the use of DFM models. Furthermore, the structure and modularity of PorePy makes it ideally suited for experimentation with mathematical models, as well as prototyping of simulation approaches.

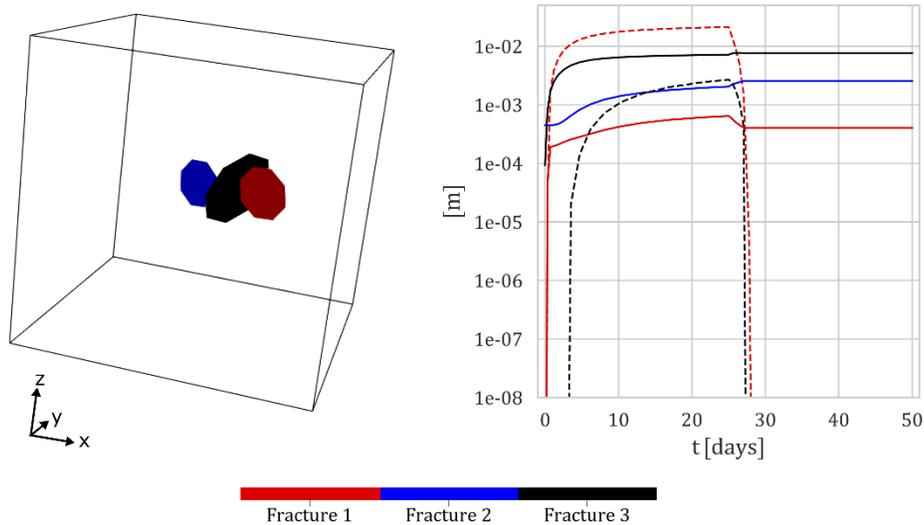

Figure 12: Left: Domain geometry with numbering of the three fractures. Fluid is injected in fracture 1 during the first 25 days, after which the well is shut. Right: $L^2$ norm normalized by fracture area of the normal (dashed lines) and tangential displacement jumps (solid lines) for each fracture.

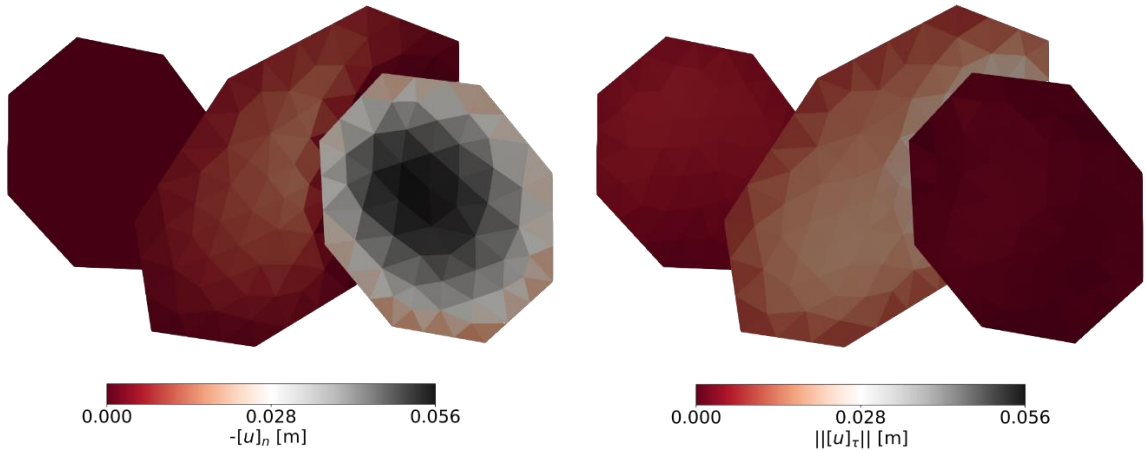

Figure 13: Normal and tangential displacements jumps on the fractures at the end of the injection phase to the left and right, respectively. The orientation of the fracture network corresponds to that in Figure 12, with the injection fracture to the right.

## 6 Conclusions

The complexity in modeling and simulation of multiphysics processes in fractured porous media, combined with a strong current research focus and corresponding developments, calls for flexible simulation tools that facilitate rapid prototyping of models and discretization methods. This paper presents design principles for simulation software for dynamics in fractured porous media, together with their implementation in the open-source simulation tool PorePy. The combined framework for modeling and simulation is based on the Discrete Fracture Matrix model, where fractures and their intersections are represented as separate lower-dimensional geometric objects. The framework facilitates flexibility for multiphysics dynamics and reuse of existing code written for non-fractured domains; hence, it is well suited for extending other software packages to mixed-dimensional problems.

The open-source software PorePy demonstrates the capabilities of the suggested framework: It provides automatic meshing of complex fracture networks in two and three dimensions, and contains implemented numerical methods for flow, transport, poroelastic deformation of the rock, and fracture deformation modeled by contact mechanics. The implementation performs well for benchmark problems in flow, poroelastic deformation and fracture deformation. Further, multiphysics simulations of fully coupled flow and non-linear transport, and of fracture deformation under poromechanical deformation of a domain demonstrates the versatility of the software.

## Acknowledgements

This work has been funded in part by Norwegian Research Council grant 250223, 244129/E20 and 267908/E20, and by a VISTA Scholarship from the Norwegian Academy of Science and Letters.